\documentclass[12pt]{article}

\usepackage{mathrsfs}
\usepackage{amssymb,amsthm,amsmath,amssymb,wrapfig,dsfont,authblk}
\usepackage{mathtools}
\usepackage{blkarray}
\usepackage{mathrsfs}
\usepackage{mathtools}
\usepackage{amsmath}
\usepackage{amsfonts}
\usepackage{amsmath,amsfonts,amssymb,amsthm,epsfig,epstopdf,titling,url,array}
\usepackage{lipsum}
\usepackage{mathtools, booktabs, bigstrut}
\usepackage{blindtext}
\usepackage{graphicx}

\usepackage{thmtools,thm-restate}
\usepackage[margin=0.8in]{geometry}

\usepackage{varioref}
\usepackage[usenames,svgnames,xcdraw,table]{xcolor}
\definecolor{DarkBlue}{rgb}{0.1,0.1,0.5}
\definecolor{DarkGreen}{rgb}{0.1,0.5,0.1}
\usepackage[backref=page]{hyperref}
\hypersetup{
     colorlinks   = true,
     linkcolor    = DarkBlue, 
     urlcolor     = DarkBlue, 
	 citecolor    = DarkGreen 
}
\renewcommand*{\backref}[1]{}
\renewcommand*{\backrefalt}[4]{%
    \ifcase #1 (Not cited.)%
    \or        (Cited on page~#2)%
    \else      (Cited on pages~#2)%
    \fi}
\usepackage[capitalise,noabbrev]{cleveref}

\usepackage{appendix}

\usepackage{algorithmic}
\usepackage{algorithm}

\usepackage{graphicx,epsfig,amsmath,latexsym,amssymb,verbatim}

\newcommand{\extra}[1]{}

\usepackage{amsmath,amssymb,amsfonts}
\usepackage{amsthm}

\newtheorem{theorem}{Theorem}

\newtheorem{proposition}[theorem]{Proposition}

\theoremstyle{remark}

\def\squareforqed{\hbox{\rlap{$\sqcap$}$\sqcup$}}
\def\qed{\ifmmode\squareforqed\else{\unskip\nobreak\hfil
\penalty50\hskip1em\null\nobreak\hfil\squareforqed
\parfillskip=0pt\finalhyphendemerits=0\endgraf}\fi}
\def\endenv{\ifmmode\;\else{\unskip\nobreak\hfil
\penalty50\hskip1em\null\nobreak\hfil\;
\parfillskip=0pt\finalhyphendemerits=0\endgraf}\fi}

\newenvironment{proof+}[1]{\noindent \textbf{{Proof #1~} }}{\qed\medskip}

\mathchardef\ordinarycolon\mathcode`\:
\mathcode`\:=\string"8000
\def\vcentcolon{\mathrel{\mathop\ordinarycolon}}
\begingroup \catcode`\:=\active
  \lowercase{\endgroup
  \let :\vcentcolon
  }

\newcommand{\nc}{\newcommand}

\nc{\barA}{\overline{A}}
\nc{\barB}{\overline{B}}
\nc{\barC}{\overline{C}}
\nc{\barD}{\overline{D}}
\nc{\barR}{\overline{R}}
\nc{\barX}{\overline{X}}
\nc{\barY}{\overline{Y}}
\nc{\barU}{\overline{U}}

\nc{\ALG}{\textsc{Alg}}

\newcommand{\A}{\mathscr{A}}
\newcommand{\B}{\mathscr{B}}

\newcommand{\X}{\mathscr{X}}
\newcommand{\E}{\mathcal{E}}
\newcommand{\M}{\mathscr{M}}
\newcommand{\LL}{\mathscr{L}}
\newcommand{\N}{\mathscr{N}}
\newcommand{\U}{\mathscr{U}}
\newcommand{\Q}{\mathcal{Q}}
\newcommand{\Y}{\mathscr{Y}}
\newcommand{\K}{\mathcal{K}}
\newcommand{\SA}{\mathscr{S\mspace{-5mu}A}}
\newcommand{\HH}{\mathscr{H}}
\newcommand{\HL}{\mathscr{H}\oplus \mathscr{L}_1}
\newcommand{\CB}{\mathcal{C}_b(\Psi)}
\newcommand{\lamZ}{\lambda(E(z))}
\newcommand{\GammaDomainRange}{\Gamma: \Omega \times \Omega\rightarrow B(\mathcal{C}_b(\Psi), B(\Y))}
\newcommand{\J}{\mathcal{J}}

\providecommand{\keywords}[1]{\textbf{\textit{Keywords---}} #1}

\begin{document}
	
	\title{{ On the Nevanlinna problem - }\\
	Characterization of all Schur-Agler class solutions affiliated with a given kernel}
	\author{Tirthankar Bhattacharyya \thanks{Department of Mathematics, Indian Institute of Science, Bangalore 560012, India. e-mail: tirtha@iisc.ac.in} \and Anindya Biswas \thanks{Department of Mathematics, Indian Institute of Science, Bangalore 560012, India. e-mail: anindyab@iisc.ac.in} \and Vikramjeet Singh Chandel \thanks{Department of Mathematics, Indian Institute of Technology Bombay, Powai, Mumbai, 400076, India. e-mail: vikramc@math.iitb.ac.in}}

	\date{}
	\maketitle

	\begin{abstract}
		\noindent
	Given a domain $\Omega$ in $\mathbb{C}^m$, and a finite set of points $z_1,z_2,\ldots, z_n\in \Omega$ and $w_1,w_2,\ldots, w_n\in \mathbb{D}$ (the open unit disc in the complex plane), the \textit{Pick interpolation problem} asks when there is a holomorphic function $f:\Omega \rightarrow \overline{\mathbb{D}}$ such that $f(z_i)=w_i,1\leq i\leq n$. Pick gave a condition on the data $\{z_i, w_i:1\leq i\leq n\}$ for such an $interpolant$ to exist if $\Omega=\mathbb{D}$. Nevanlinna characterized all possible functions $f$ that \textit{interpolate} the data. We generalize Nevanlinna's result to a domain $\Omega$ in $\mathbb{C}^m$ admitting holomorphic test functions when the function $f$ comes from the Schur-Agler class and is affiliated with a certain completely positive kernel. The Schur class is a naturally associated Banach algebra of functions with a domain. The success of the theory lies in characterizing the Schur class interpolating functions for three domains - the bidisc, the symmetrized bidisc and the annulus - which are affiliated to given kernels.
\end{abstract}

\keywords{Nevanlinna problem, Schur-Agler class, Colligation, Test function, Several complex variables}

\section{Introduction.}

Given a solvable scalar valued interpolation problem from the unit disc into the unit disc, Nevanlinna in \cite{N-29} gave a complete set of solutions. The main result of this paper, Theorem \ref{parametrization}, is a far reaching generalization of Nevanlinna's result.

\subsection{Test functions}\label{test}
A collection $\Psi$ of $\mathbb{C}$-valued functions on a set $\Omega$ is called a set of \textit{test functions} (see \cite{B_H} and \cite{D-M}) if the following conditions hold:
\begin{enumerate}
	\item $\sup\{|\psi(x)|:\psi\in \Psi\}<1$ for each $x\in \Omega$.
	\item For each finite subset $F$ of $\Omega$, the collection $\{\psi |_F:\psi\in \Psi\}$ together with the constant function generates the algebra of all $\mathbb{C}$-valued functions on $F$.
\end{enumerate}
The second condition is not essential for the development of the theory, but it makes some situations simpler (it is excluded in \cite{B_H} but not in \cite{D-M}). The collection $\Psi$ is a natural topological subspace of $\overline{\mathbb{D}}^\Omega$ equipped with the product topology. For every $x\in \Omega$, there is an element $E(x)$ in $\mathcal{C}_b(\Psi)$, {\em the $C^*$-algebra of all bounded functions on} $\Psi$, such that $E(x)(\psi)=\psi(x)$. Clearly,
\begin{math}
\| E(x)\|= \text{sup}_{\psi\in \Psi} |\psi(x)|<1
\end{math}
for each $x\in \Omega$. The functions $E(x)$ will be used at several places in this paper.

\subsection{Completely Positive Kernels}\label{CPK} 
A positive kernel $k$ on a set $\Omega$ is a function $k:\Omega\times \Omega \rightarrow \mathbb{C}$ such that for any $n\geq 1$, any $n$ points $z_1,z_2,\ldots,z_n$ in $\Omega$ and any $n$ complex numbers $c_1,c_2,\ldots, c_n$, we have
\begin{align*}
\sum_{i=1}^{n}\sum_{j=1}^{n}\overline{c_i}c_j k(z_i,z_j)\geq 0.
\end{align*}
If $\mathcal{E}$ is a Hilbert space and $k:\Omega\times \Omega\rightarrow B(\mathcal{E})$ is a function, then $k$ is called a positive kernel if for any $n\geq 1$, any $n$ points $z_1,z_2,\ldots,z_n$ in $\Omega$ and any $n$ vectors $\mathbf{e_1,e_2,\ldots, e_n}$ in $\mathcal{E}$, we have
\begin{align}\label{op. valued kernel}
\sum_{i=1}^{n}\sum_{j=1}^{n}\langle k(z_i,z_j)\mathbf{e_j}, \mathbf{e_i} \rangle\geq 0.
\end{align}
The concept of a positive kernel does not cease here. Let $\A$ and $\B$ be two $C^*$-algebras and let $\Gamma$ be a function on $\Omega \times \Omega$ taking values in $B(\A, \B)$ (space of all bounded linear operators from $\A$ to $\B$). $\Gamma$ is called a \textit{completely positive kernel} if 
\begin{align} \label{compositivekernel}
	\sum_{i,j =1}^{n} b_i^* \  \Gamma(z_i, z_j)(a_i ^*a_j) b_j \geq 0
\end{align}
for all $n \geq 1$, $a_1, a_2, \dots, a_n \in \A$, $b_1, b_2, \dots, b_n \in \B$ and $z_1, z_2, \dots, z_n \in \Omega$.\\

It is well-known that a completely positive kernel has a Kolmogorov decomposition (see \cite{B-B-F-t} or \cite{B_H} or \cite{B-B-L-S}). The following situation will occur several times in our paper. Let us consider a Hilbert space $\Y$ and a collection $\Psi$ of test functions. Take the  $C^*$-algebras $\mathcal{C}_b(\Psi)$ and $B(\Y)$, a Hilbert space $\X$, a unital $*$-representation $\mu :\CB \rightarrow B(\X)$ and a function $h:\Omega\rightarrow B(\X,\Y)$. Then it is easy to see that
$$\Gamma(x,y)(\delta)=h(z)\mu(\delta)h(w)^*\,\,\text{for all}\,\, z,w\in\Omega\,\,\text{and}\,\, \delta\in \CB$$
is a completely positive kernel. Conversely, if
 $$\Gamma: \Omega \times \Omega\rightarrow B(\mathcal{C}_b(\Psi), B(\Y))$$
  is a completely positive kernel, then there exists a Hilbert space $\X$, a unital $*$-representation $\mu: \mathcal{C}_b(\Psi) \rightarrow B(\X)$ and a function $h: \Omega \rightarrow B(\X, \Y)$ such that
\begin{align}\label{Kolmogorv Decomposition}
\Gamma(z, w)(\delta) = h(z) \mu(\delta) h(w)^* \,\,\text{for all} \ \ z,w \in \Omega \ \ \text{and} \ \ \delta \in \mathcal{C}_b(\Psi).
\end{align}
This is called the \textit{Kolmogorov decomposition.}

For example, take a Hilbert space $\Y$ and a function $\psi:\Omega \rightarrow \mathbb{C}$ from a given collection $\Psi$ of test functions. Since $|\psi(z)|<1$ for each $z\in \Omega$, for any finite subset $\Omega_0\subset \Omega$ we have $sup_{z\in\Omega_0} |\psi(z)|<1$. Now take the map $\Gamma_{\psi}:\Omega \times \Omega \rightarrow B(\CB, B(\Y))$ defined by
\begin{align}\label{Gamma_Psi}
\Gamma_{\psi} (z,w) (\delta)= \frac{\delta(\psi)}{1- \psi(z) \overline{\psi(w)}} I_\Y,\,\, z,w\in\Omega,\,\,\delta \in \CB.
\end{align}
Since the positivity condition (\ref{compositivekernel}) involves only a finite number of points from $\Omega$, we shall have no difficulty in checking the positivity of this map. Let $\X =\bigoplus_{j=0} ^{\infty}\Y_j$ where $\Y_j =\Y$ for every $j$. Define $\mu_\psi :\CB \rightarrow B(\X)$ by
$$\mu_\psi (\delta)= \delta (\psi)I_{\X}$$
and $h_\psi :\Omega \rightarrow B(\X, \Y)$ by
$$h_\psi (z)\big(\oplus_{j=0} ^{\infty} y_j\big)= \sum_{j=0}^{\infty} \psi(z)^j y_j.$$
Clearly $\mu_\psi$ is a unital $*$-representation and $h_\psi (z)$ is bounded linear transformation. Note that
$$\Gamma_{\psi} (z,w)(\delta) = h_\psi (z)\mu_\psi(\delta)h_\psi (w)^*.$$
Hence by (\ref{Kolmogorv Decomposition}), $\Gamma_{\psi}$ is a completely positive kernel.

\subsection{$\Psi$-unitary Colligations}
Let $\X$, $\U$ and $\Y$ be Hilbert spaces and let $\Psi$ be a fixed set of test functions. By a $\Psi$-unitary colligation, we mean a pair $(U, \rho)$ where $U$ is a unitary operator from $\X \oplus \U$ to $\X \oplus \Y$, and $ \rho: \mathcal{C}_b(\Psi) \rightarrow B(\X)$ is a $*$-representation. If we write $U$ as
\[
U =
\bordermatrix{ & \X & \U \cr
	\X & A & B \cr
	\Y & C & D}, \qquad
\]
then we can define a bounded $B(\U, \Y)$ valued function on $\Omega$, given by
\begin{align}
  f(x) = D + C \rho (E(x)) (I_{\X}- A \rho(E(x)))^{-1} B \ \ \forall \ x \in \Omega,
\end{align}
equivalently,
\begin{align}
  f(x) = D + C(I_{\X}- \rho (E(x)) A)^{-1} \rho (E(x)) B  \ \ \forall \ x \in \Omega.
\end{align}
This $f$ is called the transfer function associated with $(U,\rho)$. Since $U^*$ is also a unitary, we have that
\begin{align*}
g(x) = D^* + B^*(I_{\X}- \rho (E(x)) A^*)^{-1} \rho(E(x)) C^*
\end{align*}
 is the transfer function of the colligation $(U^*,\rho)$.

\subsection{The $\Psi$-Schur-Agler Class}\label{SAClass}
Let $\E$ be a Hilbert space and $\Omega$ an abstract set. We consider a $B(\E)$-valued kernel $K$ (satisfying (\ref{op. valued kernel})) on $\Omega$. For this $K$, there is a Hilbert space $\mathcal{H}(K)$ of $\E$-valued functions on $\Omega$ such that span of the set
\begin{align*}
\{K(\cdot, \omega)\mathbf{e}:\mathbf{e}\in \E,\, \omega\in \Omega \}
\end{align*}
is dense in $\mathcal{H}(K)$ and for any $\mathbf{e}\in \E$, $\omega \in \Omega$ and $h\in \mathcal{H}(K)$, we have
\begin{align*}
\langle h, K(\cdot, \omega)\mathbf{e}\rangle_{\mathcal{H}(K)} = \langle h(\omega), \mathbf{e} \rangle_{\E}.
\end{align*}

Given a set of test functions $\Psi$ on $\Omega$, a kernel $K:\Omega \times \Omega\rightarrow B(\E)$ is said to be {\em $\Psi$-admissible} if the map $M_\psi$, sending each element $h\in \mathcal{H}(K)$ to $\psi\cdot h$, is a contraction on $\mathcal{H}(K)$. We denote the set of all $B(\E)$-valued $\Psi$-admissible kernels by $\K_\Psi(\E)$. For two Hilbert spaces $\U$ and $\Y$, we say that {\em $S: \Omega \rightarrow B(\U, \Y)$ is in $H^\infty_\Psi(\U, \Y)$} if there is a non-negative constant $C$ such that the $B(\Y \otimes \Y)$-valued function
\begin{equation} \label{C} (C^2 I_\Y - S(x)S(y)^*)\otimes k(x,y)  \end{equation}
 is a positive $B(\Y\otimes \Y)$-valued kernel for every $k$ in $\K_\Psi(\Y)$. If $S$ is in $H^\infty_\Psi(\U, \Y)$, then we denote by $C_S$ the smallest $C$ which satisfies \eqref{C}. The collection of maps $S\in H^\infty_\Psi(\U, \Y)$ for which $C_S$ is no larger than $ 1$ is called the $\Psi$-Schur-Agler class and it is denoted by $\SA_\Psi (\U,\Y)$.\\

The plan of the paper is as follows. In Section \ref{Char-schur-agler}, we prove a characterization of functions in the class $\SA_\Psi (\U,\Y)$. This is followed by Section \ref{Holomorphic test functions} which describes the usefulness of taking holomorphic test functions. Section \ref{Affiliated} consists of the description of  an auxiliary function $G$. Section \ref{Main Result} has the main theorem. In Section \ref{Examples}, we give applications of our main result.

 \section{Characterization of $\SA_\Psi (\U,\Y)$}\label{Char-schur-agler}

Variants of the following theorem exist in various forms in literature, see \cite{D-M} and \cite{B_H} and the references therein. We did not find it in the form that we shall need. The most non-trivial implication is $\mathbf{1}. \Rightarrow \mathbf{2}.$ and we shall prove this since we did not find, in the literature, a proof of it. Other implications are easy to see.

\begin{theorem}\label{RealTheo}
	Consider a function $S_0$ on some subset $\Omega_0$ of $\Omega$ with values in $B(\U, \Y)$. Then the following conditions are equivalent.
	\begin{enumerate}
		\item There exists an $S$ in $\SA_\Psi(\U,\Y)$ such that $S|_{\Omega_0} =S_0$.
		\item $S_0$ has an {\em Agler decomposition} on $\Omega_0$, that is, there exists a completely positive kernel $\Gamma: \Omega_0 \times \Omega_0 \rightarrow B(\mathcal{C}_b(\Psi), B(\Y))$ so that
		\begin{align*}
		I_{\Y} - S_0(z) S_0(w)^* = \Gamma(z,w) (1- E(z) E(w) ^*) \ \ \text{for all} \ z, w \in \Omega_0.
		\end{align*}
		\item There exists a Hilbert space $\X$, a $*$-representation $\rho: \mathcal{C}_b(\Psi) \rightarrow B(\X)$ and a $\Psi$-unitary colligation $(V, \rho)$ such that writing $V$ as
		\[
		V =
		\bordermatrix{ & \X & \U \cr
			\X & A & B \cr
			\Y & C & D}, \qquad
		\]
		one has
		\begin{align}
		S_0(z) = D+ C(I_{\X}- \rho(E(z)) A) ^{-1} \rho(E(z)) B \ \ \text{for all} \ z \in \Omega_0.
		\end{align}
		\item There exists a Hilbert space $\X$, a $*$-representation $\rho: \mathcal{C}_b(\Psi) \rightarrow B(\X)$ and a $\Psi$-unitary colligation $(W, \rho)$ such that writing $W$ as
		\[
		W =
		\bordermatrix{ & \X & \Y \cr
			\X & A_1 & B_1 \cr
			\U & C_1 & D_1}, \qquad
		\]
		one has
		\begin{align}\label{transfer}
		S_0(z)^* = D_1+ C_1 (I_{\X}- \rho(E(z))^* A_1) ^{-1} \rho(E(z))^* B_1 \ \ \text{for all} \ z \in \Omega_0.
		\end{align}
	\end{enumerate}
\end{theorem}

\textbf{Proof.} As mentioned before, we shall only prove that $\mathbf{1}.$ implies $\mathbf{2}.$
Consider an $S\in \SA_\Psi(\U, \Y)$ and a $\Psi$-admissible kernel $K:\Omega \times \Omega\rightarrow B(\Y)$. As is usual, denote
\begin{align*}
\mathcal{H}(K)=\overline{span}\{K(\cdot, w)y:w\in \Omega, y\in \Y\}.
\end{align*}
Define a linear transformation $T^*$ on the dense subspace $span \{K(\cdot, w)y_1 \otimes y_2:w\in \Omega, y_1, y_2\in \Y\}$ by first defining $$T^*\Big (K(\cdot, w)y_1 \otimes y_2\Big) =K(\cdot, w)y_1 \otimes S(w)^* y_2,$$
and then extending linearly.
For $w_i\in \Omega, \,y_{1i}, y_{2,i}\in \Y, 1\leq i\leq n,$ we have
\begin{align*}
&||\sum_{i=1}^{n}K(\cdot, w_i)y_{1i} \otimes y_{2i}||^2-||T^* \Big(\sum_{i=1}^{n}K(\cdot, w_i)y_{1i} \otimes y_{2i}\Big)||^2\\ &=\sum_{i,j=1}^{n}\Bigg \langle \Big( (I_\Y -S(w_j)S(w_i)^*)\otimes K(w_j,w_i)\Big ) y_{1i} \otimes y_{2i}, y_{1j} \otimes y_{2j}\Bigg \rangle
\end{align*}
Since $S\in \SA_\Psi(\U, \Y)$, the last expression is nonnegative. So $T$ defines a contraction from $\mathcal{H}(K)\otimes \U$ to  $\mathcal{H}(K)\otimes \Y.$ We need the following lemma to continue with the proof.

\textbf{Lemma.} Let $J : \Omega \times \Omega\rightarrow B(\Y)$ be a self-adjoint function, i.e., $J(z,w)= J(w,z)^*$ for all $z,w\in \Omega$ (see page 174 of \cite{A-M}). If
\begin{align}\label{TheJ}
J\oslash K: (z,w)\mapsto J(z,w)\otimes K(z,w)
\end{align}
is a positive kernel for every $B(\Y)$-valued admissible kernel $K$, then there is a completely positive kernel $\GammaDomainRange$ such that
$$J(z,w)=\Gamma(z,w)(1- E(z)E(w)^*)\,\, \text{for all}\,\, z,w\in \Omega.$$
\textbf{Proof of the lemma.}\\
We prove the result for a finite subset $\Omega_0=\{w_1,w_2,\ldots, w_n\}$ of $\Omega$ and apply Kurosh's theorem.\\
Consider the following subset of $n\times n$ self-adjoint operator matrices with entries in $B(\Y)$
\begin{align*}
W_{\Omega_0}=&\Big\{\Big(\Gamma(w_i,w_j)(1-E(w_i)E(w_j)^*)\Big)_{1\leq i,j\leq n}:\\ &\Gamma:\Omega_0\times \Omega_0\rightarrow B(\CB, B(\Y))
\text{ is a completely positive kernel}\Big\}
\end{align*}
Clearly $W_{\Omega_0}\subset B(\Y^n)$ is a convex set and it is invariant under multiplication by positive real scalars in the space of $n\times n$ self-adjoint matrices. Such a set is called a \textit{wedge} (see \cite{A-M}, page 169)).

It needs to be shown that  $W_{\Omega_0}$ is closed in the $weak^*$-topology of $ B(\Y^n).$ To that end, start with a net $\Big[ \Gamma_\nu (w_i,w_j)(1-E(w_i)E(w_j)^*)\Big]_{1\leq i,j\leq n} = \Big[ A_\nu \Big]_{1\leq i,j\leq n}$ in $W_{\Omega_0}$ and suppose that it converges in the $weak^*$-topology to an $n\times n$ self-adjoint matrix $A=\Big[A_{ij}\Big]$ with entries in $B(\Y)$. This means that for every $X=\Big[X_{lp} \Big]$ in $B_1 (\Y^n)$ (the space of trace class operators on $\Y^n$), $\{tr(A_\nu X)\}$ converges to $tr (AX)$. Let $u,v\in \Y$ with $||u||,||v||\leq 1$ and choose $X$ to be the operator matrix which has $u\otimes v$ as its $(j,i)$-th entry and zeros elsewhere. Then $tr(A_\nu X)=\Big\langle A_{\nu}u,v \Big\rangle$ tends to $tr(A X)=\Big\langle Au,v \Big\rangle$. We have that
\begin{align*}
\Big\langle\Gamma_\nu (w_i,w_j)(1-E(w_i)E(w_j)^*)u,v \Big\rangle\,\, \longrightarrow\,\, \Big\langle A_{ij} u,v \Big\rangle.
\end{align*}
Now $1- E(w_i)E(w_i)^*\geq 1-||E(w_i)||^2> 0$ gives us that there is an $\epsilon>0$ such that $1-E(w_i)E(w_i)^*> \epsilon\cdot 1\,\, \text{for all}\,\, i=1,2,\ldots, n.$ Hence we get that
\begin{align*}
\Big\langle\Gamma_\nu (w_i,w_i)(1-E(w_i)E(w_i)^*)u,u \Big\rangle \geq \epsilon\, \Big\langle\Gamma_\nu (w_i,w_j)u,u \Big\rangle \,\, \text{for all}\,\, i=1,2,\ldots, n.
\end{align*}
Since the left hand side converges, we can find an $M>0$ such that
\begin{align*}
sup_{\nu}\,\, \Big\langle\Gamma_\nu (w_i,w_i)u,u \Big\rangle \leq M \,\, \text{for all}\,\, i=1,2,\ldots, n.
\end{align*}
Also for any $\delta \in \CB,$ we have
\begin{align*}
\Big\langle\Gamma_\nu (w_i,w_i)(\delta \delta^*)u,u  \Big\rangle \leq ||\delta||^2 \Big\langle\Gamma_\nu (w_i,w_i)u,u \Big\rangle \leq M ||\delta||^2.
\end{align*}
Now we need the following result which is a Cauchy-Schwarz type inequality. Since we could not find an exact reference, we are giving a proof of the claim.

\textbf{Claim:} If $\Gamma$ is completely positive, $\delta \in \CB$, $z,w \in \Omega$ and $u,v\in \Y$, then we have
\begin{align*}
\Big|\Big\langle\Gamma (z,w)(\delta \delta ^*)u,v \Big\rangle\Big|^2 \leq \Big\langle\Gamma (z,z)(\delta \delta ^*)v,v \Big\rangle   \Big\langle\Gamma (w,w)(\delta \delta ^*)u,u \Big\rangle.
\end{align*}
\textbf{Proof of the claim:}\\
Note that, using Kolmogorov decomposition (\ref{Kolmogorv Decomposition}), we can find a Hilbert space $\X_1$, a unital $*$-representation $\mu: \mathcal{C}_b(\Psi) \rightarrow B(\X_1)$ and a function $h: \Omega \rightarrow B(\X_1, \Y)$ such that
\begin{align*}
\Gamma(z, w)(\delta) = h(z) \mu(\delta) h(w)^* \,\,\text{for all} \ \ z,w \in \Omega \ \ \text{and} \ \ \delta \in \mathcal{C}_b(\Psi).
\end{align*}
So we have
\begin{align*}
\Big|\Big\langle\Gamma (z,w)(\delta \delta ^*)u,v \Big\rangle\Big|^2 &=\Big|\Big\langle h(z)\mu(\delta \delta ^*)h(z)^*u,v \Big\rangle\Big|^2 \\
&=\Big|\Big\langle \mu(\delta ^*)h(z)^*u,\mu(\delta ^*)h(w)^*v \Big\rangle\Big|^2\\
&\leq \Big\langle \mu(\delta ^*)h(z)^*u,\mu(\delta ^*)h(z)^*u \Big\rangle \Big\langle \mu(\delta ^*)h(w)^*v,\mu(\delta ^*)h(w)^*v \Big\rangle.
\end{align*}
The last expression equals to
$$\Big\langle\Gamma (z,z)(\delta \delta ^*)v,v \Big\rangle   \Big\langle\Gamma (w,w)(\delta \delta ^*)u,u \Big\rangle.$$
 Hence the claim follows.\qed
Using this result, we conclude that
\begin{align*}
\Big|\Big\langle\Gamma_\nu (w_i,w_j)(\delta \delta^*)u,v  \Big\rangle\Big|\leq M ||\delta||^2\,\,\text{for every}\,\, i,j=1,2,\ldots, n.
\end{align*}
Therefore, for each $\delta\in \CB$, $u,v\in \Y$ and $i,j=1,2,\ldots, n$, the net $\Big\{\Big\langle\Gamma_\nu (w_i,w_j)(\delta \delta^*)u,v  \Big\rangle\Big\}$
is bounded. Since $\Omega_0$ is finite, we get a subnet ${\nu_l}$ such that $\Big\{\Big\langle\Gamma_{\nu_l} (w_i,w_j)(\delta \delta^*)u,v  \Big\rangle\Big\}$
converges to some number depending on $\delta, u$ and $v$, $\text{for every}\,\, i,j=1,2,\ldots, n$. Define a completely positive kernel
$ \Gamma :\Omega_0 \times \Omega_0 \rightarrow B(\CB, B(\Y))$
by $\Big\langle\Gamma(w_i,w_j)(\delta )u,v  \Big\rangle=\lim_l\Big\langle\Gamma_{\nu_l} (w_i,w_j)(\delta)u,v  \Big\rangle$ and extend it trivially to the whole set $\Omega \times \Omega.$ Consequently, for every $u,v\in \Y$ we have
\begin{align*}
\Big\langle\Gamma(w_i,w_j)(1-E(w_i)E(w_j)^* )u,v  \Big\rangle=\Big \langle A_{ij}u, v \Big \rangle.
\end{align*}
This proves that $W_{\Omega_0}$ is $weak*$-closed.

The $n\times n$  matrix $\Big [I_\Y\Big]$ with each entry equal to $I_\Y$, is in $W_{\Omega_0}$. Indeed, let $\psi \in \Psi$ and take $\Gamma_{\psi}:\Omega_0 \times \Omega_0 \rightarrow B(\CB, B(\Y))$ defined by
\begin{align*}
\Gamma_{\psi} (z,w) (\delta)= \frac{\delta(\psi)}{1- \psi(z) \overline{\psi(w)}} I_\Y,\,\, z,w\in\Omega_0.
\end{align*}
Note that $\sup_{z\in \Omega_0} |\psi(z)|<1$ as $\Omega_0$ is finite. From (\ref{Gamma_Psi}), it is clear that $\Gamma_{\psi}$ is completely positive.
So we have
\begin{align*}
 \Gamma_{\psi} (z,w)(1-E(z)E(w)^*)=I_\Y
\end{align*}
and hence we conclude that $\Big [I_\Y\Big]\in W_{\Omega_0}$.

Also, the restriction of $J$  (see (\ref{TheJ})) on $\Omega_0 \times \Omega_0$, that is, $J|_{\Omega_0 \times \Omega_0}=\mathcal{J}$ is in $W_{\Omega_0}$. If possible, let $\J \notin W_{\Omega_0}$. By Theorem $3.4$ in \cite{Rudin}, we get a $weak^*$-continuous linear functional $L$ on $B(\Y^n)$ whose real part is nonnegative on $W_{\Omega_0}$ and strictly negative at $\J$. We replace $L(R)$ by $\frac{L(R) + \overline{L(R)}}{2},\,\, R\in B(\Y^n), $ and denote it by $L$ itself. Since $L$ is $weak^*$-continuous and for any locally convex space $X$, we have $(X^*; weak^*)^*=X$ (see Theorem V.1.3 in \cite{Conway}), we find that $L$ is of the form $L(R) = tr(RC)$ for some $n \times n$ self-adjoint compact $C\in B_1 (\Y^n)$ whose entries are in the ideal of trace class operators on $Y^n$.

Let $\{e_n: n \geq 1\}$ be an orthonormal basis of $\Y$. Given a bounded operator $A$ on $Y$, define its $transpose$ to be the linear transformation on $\Y$ whose matrix entries with respect to the basis above are  $\langle A^t e_j, e_i \rangle = \langle A e_i, e_j \rangle$. It is easy to see that this defines a bounded operator. Indeed, if $u \in \Y$ is given by $u = \sum u_i e_i,$ then we define $\overline{u} = \sum \overline{u_i} e_i$ and then we have
$ \langle A^t u, v  \rangle =  \langle A \overline{v}, \overline{u}  \rangle $ for $u, v \in Y$ which on application of the Cauchy-Schwarz inequality yields boundedness.

Since $C$ obtained above is a block $n \times n$ operator matrix $C = \left( \left( C(w_i, w_j) \right) \right)_{i,j=1}^n,$ define $C^t$ to be the block $n \times n$ operator matrix whose $(i,j)$th. entry is $ C^t (w_i, w_j) = C (w_j, w_i)^t $. In other words, $ \langle C^t (w_i, w_j) u, v  \rangle = \langle C (w_j, w_i)\overline{v}, \overline{u} \rangle $ for $u, v \in Y$.

We shall show that $C^t$ is a $B(\Y)$-valued positive kernel on $\Omega_0$, that is, for $u_i \in \Y, 1\leq i\leq n$,
$$\sum_{i,j =1}^{n} \big \langle C^t (w_i, w_j) u_j, u_i \big \rangle \geq 0 $$

To see this, note that
$$ \sum_{i,j =1}^{n} \big \langle C^t (w_i, w_j) u_j, u_i \big \rangle =   \sum_{i,j =1}^{n} \big \langle \overline{u_i},  C (w_j, w_i)^* \overline{u_j} \big \rangle =   \sum_{i,j =1}^{n} \big \langle tr( (\overline{u_i} \otimes \overline{u_j}) C_{j i}  ) \big \rangle. $$
This last quantity is $tr(D C)$, where $D = (D_{i,j}) = (\overline{u_i} \otimes \overline{u_j})$, and hence equals $L(D)$.

For any function $u : \Omega \rightarrow \Y$, if $\Delta_\psi : \Omega \times \Omega \rightarrow B ( \mathcal{C}_b(\Psi), \mathcal{B}(\Y) )$ is the function defined by
$$ \Delta_\psi (z,w) : \delta \mapsto \frac{\delta(\psi)  \ u(z) \otimes u(w) }{1 - \psi(z) \overline{\psi(w)}}  ,\,\, z,w\in \Omega,$$
then
$$ \Delta_{\psi} (z,w) (\delta_1 \overline{\delta_2})= \frac{ \delta_1 (\psi) \overline{\delta_2 (\psi)}}{1 - \psi(z) \overline{\psi(w)} } \ u(z) \otimes u(w).  $$
Let $\delta_1,\delta_2,\ldots, \delta_n \in \mathcal{C}_b(\psi)$ , $B_1, B_2,\ldots,B_n \in B(\Y)$, $v \in \Y$ and $z_1,z_2, \ldots, z_n \in \Omega$. Then we have
$$   \bigg \langle \Big(  \sum_{i,j =1}^{n} B^* _i \Delta_\psi (z_i, z_j)  (\delta_i^* \delta_j \big)B_j\Big) v, v \bigg \rangle = \sum_{i,j =1}^{n} \frac{ \alpha_i \overline{\alpha_j}}{1 - \psi(z_i) \overline{ \psi(z_j)}} $$
where $\alpha_i = \overline{\delta_i(\psi)} \big \langle u(z_i), B_i(v) \big \rangle  \ \ i = 1,2, \ldots, n.$
The last expression is clearly non-negative. Therefore, $\Delta_{\Psi}$ is completely positive.
Now, $\Delta_{\psi} ( 1 - E(z) E(w)^* ) = u(z) \otimes u(w) $ gives that $D$ is in $W_{\Omega_0}$.
Hence,  $ \sum_{i,j =1}^{n} \big \langle C^t (w_i, w_j) u_j, u_i \big \rangle = L(D) \geq 0.$ Thus, $C^t$ is a $B(\Y)$-valued positive kernel on $\Omega_0$.

In fact, $C^t$ is admissible. To see that, consider the function $\Xi_\psi :\Omega \times \Omega\rightarrow B(\CB, B(\Y))$ defined by
\begin{align*}
\Xi_\psi (z,w)(\delta)= \delta(\psi)u(z)\otimes u(w)
\end{align*}
for $\psi \in \Psi$ and $u:\Omega \rightarrow \Y$. Let $\delta_i\in \CB,B_i\in \B(\Y)$ and $z_i\in \Omega, 1\leq i\leq n.$\\
Then for any $v\in \Y$ we have
\begin{align*}
\sum_{i,j = 1}^{n}\Big \langle \Big(B_i ^* \Xi_\psi (z_i,z_j)(\delta_i ^* \delta_j)B_j \Big)v,v\Big \rangle = |\sum_{i= 1}^{n}\delta_i(\psi)\langle B_iv, u(z_i)\rangle|^2 \geq 0.
\end{align*}
So $\Xi_\psi$ is completely positive and hence $$\Xi_\psi (z,w)(1-E(z)E(w)^*)|_{\Omega_0 \times \Omega_0} \in W_{\Omega_0}.$$
Now let $u_1,u_2,\ldots, u_n\in \Y$. A little computation gives
\begin{align*}
\sum_{i,j =1}^{n} \big \langle (1 - \psi(w_i) \overline{ \psi(w_j)})C^t (w_i, w_j) u_j, u_i \big \rangle =tr(\Xi_\psi C).
\end{align*}
But then $\Xi_\psi (z,w)(1-E(z)E(w)^*)|_{\Omega_0 \times \Omega_0} \in W_{\Omega_0}$ gives us $tr(\Xi_\psi C)\geq 0.$ So $C^t$ is $\Psi$-admissible.

By our assumption, the $B(\Y \otimes \Y)$-valued function $\J \oslash C^t$  on $\Omega_0 \times \Omega_0$ is positive. So for any $u_i\in \Y,1\leq i\leq n,$ we have
\begin{align*}
\sum_{i,j = 1}^{n} \Big\langle \J \oslash C^t (w_i, w_j)u_j, u_i\Big\rangle \geq 0.
\end{align*}
With $\{e_p: p \geq 1\}$ an orthonormal basis of $\Y$, set $u_i = \sum_{m=1}^{N} e_m \otimes e_m,\,\, \text{for all}\,\, i,$ where $N$ is some fixed natural number. So we have
\begin{align*}
\sum_{i,j = 1}^{n} \Big\langle \J \oslash C^t (w_i, w_j)u_j, u_i\Big\rangle
= \sum_{i,j = 1}^{n} \sum_{p,q =1}^{N} \Big\langle \J (w_i, w_j)e_p, e_q\Big\rangle \Big\langle  C^t (w_i, w_j)e_p, e_q\Big\rangle
\end{align*}
and this is nonnegative. This holds for any $N\geq 1$. On the other hand
$$ L(\J) = \sum_{i,j = 1}^{n} tr\Big(\J(w_i, w_j) C(w_j,w_i)\Big) =  \sum_{i,j = 1}^{n} \sum_{p,q =1}^{\infty} \Big\langle \J (w_i, w_j)e_p, e_q\Big\rangle \Big\langle  C^t (w_i, w_j)e_p, e_q\Big\rangle $$
which is nonnegative by the previous argument. But this is a contradiction since $L(\J)<0$. Hence $\J \in W_{\Omega_0}$.

It is easy to see that conditions in Kurosh's theorem are satisfied (see Theorem 2.56 in \cite{A-M}, page 74-75 in \cite{Arkhangel}). So the finiteness condition on the set $\Omega_0$ can be removed.

This completes the proof of the statement that under given conditions there is a completely positive kernel $\GammaDomainRange$ such that
$$J(z,w)=\Gamma(z,w)(1- E(z)E(w)^*)\,\, \text{for all}\,\, z,w\in \Omega.$$
This completes the proof of the lemma.\qed

 To complete the proof of the theorem, note that if $S\in \SA_\Psi (\U,\Y)$, then clearly $(I_\Y -S(z)S(w)^*)\otimes K(z,w)$ is positive for every $\Psi$-admissible kernel $K$. Hence an application of the result above shows that there is a completely positive kernel $\GammaDomainRange$ such that
\begin{align*}
I_\Y -S(z)S(w)^* = \Gamma(z,w)(1- E(z)E(w)^*)\,\,\, \text{for all}\,\,\, z,w\in \Omega.
\end{align*}

\section{Holomorphic test functions}\label{Holomorphic test functions}

When test functions are holomorphic, strong consequences can be shown to follow. That is the purpose of this section.

There are some known cases like the bidisc (\cite{A-M}), the symmetrized bidisc (\cite{B-S}) and the annulus (\cite{D-M}) where the collection of test functions has the nice property of being holomorphic. Moreover, in all these cases, we can find a point in the domain at which each all test functions vanish. We now want our collection $\Psi$ to have this property without changing the $\Psi$-Schur-Agler class. Our purpose will be served if we can find another collection $\Theta$ of holomorphic test functions having a common zero such that
$$\mathscr{S\mspace{-5mu}A}_\Psi  (\U, \Y)= \mathscr{S\mspace{-5mu}A}_\Theta  (\U, \Y).$$
Now, by definition $\mathscr{S\mspace{-5mu}A}_\Psi  (\U, \Y)$ is the set $$\{S: \Omega \rightarrow B(\U, \Y): ( I_\Y - S(x)S(y)^*)\otimes k(x,y)\,\,\text{is positive for every}\,\, k\in \K_\Psi(\Y)\}.$$
So if we can show that for a collection $\Theta$ of holomorphic test functions (with our required property) the equality $\K_\Psi(\Y)= \K_\Theta (\Y)$ holds, then we are done. The following proposition provides us with such a collection.
\begin{proposition}\label{Modified Test Functions}
	Let $\Psi$ be a collection of holomorphic test functions on a domain $\Omega \subset \mathbb{C}^m$, $w_0\in\Omega$ a point and $\Y$ a Hilbert space. Then there is a collection $\Theta=\{\theta_\psi: \psi \in \Psi\}$ such that the following conditions are met.
	\begin{enumerate}
		\item $\theta_\psi(w_0)=0$ for all $\psi \in \Psi$.
		\item $\Theta$ is a collection of holomorphic test functions.
		\item $\K_\Psi(\Y)= \K_\Theta (\Y)$
	\end{enumerate}
\end{proposition}
\textbf{Proof.} Fix a $w_0 \in \Omega$ and for each $\psi \in \Psi$ define $\varphi_\psi : \mathbb{D}\rightarrow \mathbb{D}$ by
$$\varphi_\psi (\tilde{z})=\frac{\psi(w_0)- \tilde{z}}{1- \overline{ \psi(w_0)}\tilde{z}}.$$
Since $||E(w_0)||<1$, $\varphi_\psi$ is non-constant. Let us define a holomorphic function $\theta_\psi : \Omega\rightarrow \mathbb{D}$ by $\theta_\psi = \varphi_\psi \circ \psi$. Clearly $\theta_\psi(w_0)=0$ for all $\psi \in \Psi$.

We shall show that $\Theta=\{\theta_\psi: \psi \in \Psi\}$ is a collection of holomorphic test functions. To prove the first defining condition of the test functions, we start with the assumption that for some $z_0\in \Omega$, $\sup_{\psi \in \Psi} |\theta_\psi (z_0)|=1$. So there is a sequence $\{\psi_n\}\subset \Psi$ such that $\lim_{n\rightarrow \infty}| \varphi_{\psi_n}\circ \psi_n(z_0)|=1.$ Let $f_n : \Omega \rightarrow \mathbb{D}$ denote the function $\varphi_{\psi_n}\circ \psi_n$. Since $\{f_n\}$ is a uniformly bounded collection of holomorphic functions on $\Omega$, by Montel's theorem, there is a subsequence $\{f_{n_l}\}$ of $\{f_n\}$ and a holomorphic function $f:\Omega\rightarrow \overline{\mathbb{D}}$ such that $f_{n_l}\rightarrow f$ uniformly on every compact subset of $\Omega$ as $l\rightarrow \infty$. Clearly $|f(z_0)|=1$. Since $\Omega$ is a domain (open connected set), by maximum modulus theorem, $f$ must be constant on $\Omega$. But we also have
\begin{align*}
f(w_0)= \lim_{l\rightarrow \infty} f_{n_l}(w_0)
=\lim_{l\rightarrow \infty} \varphi_{\psi_{n_l}}\circ \psi_{n_{l}}(w_0)
=0.
\end{align*}
This is a contradiction. Hence $\sup_{\psi \in \Psi} |\theta_\psi (z)|<1$ for each $z\in \Omega$.

Let $F$ be a finite subset of $\Omega$. Choose $\psi\in \Psi$ and consider $\theta_\psi = \varphi_\psi \circ \psi$. Since $F$ is finite and for each $z\in \Omega$, $\sup_{\psi\in \Psi} |\psi (z)|<1$ and $\sup_{\psi\in \Psi} |\theta_\psi (z)|<1$, there is an $\epsilon\in (0,1)$ such that $\sup_{F} |\psi| <\epsilon$ and $\sup_{F} |\theta_\psi| <\epsilon$. Now, $\varphi_\psi$ has a power series representation which converges absolutely and uniformly in $\overline{B(0,\epsilon)}\subset \mathbb{D}$. So, $\theta_\psi|_F$ is in the closed algebra generated by $\{1,\psi|_F\}$. Since $\psi =\varphi_\psi ^{-1} \circ \theta_\psi$, a similar argument gives us that $\psi|_F$ is in the closed algebra generated by $\{1,\theta_\psi|_F\}$. Hence the closed algebras generated by $\{1, \psi|_F: \psi \in \Psi\}$ and $\{1, \theta_\psi |_F: \psi\in \Psi\}$ coincide. Since $\Psi$ is a collection of test functions, the closed algebra generated by $\{1, \theta_\psi |_F: \psi\in \Psi\}$ is the algebra of all $\mathbb{C}$-valued functions on $F$. This proves that $\Theta$ is a collection of holomorphic test functions.

Now consider the Hilbert space $\Y$. $\K_\Psi(\Y)$ and  $\K_\Theta (\Y)$ are the sets of all $\Psi$-admissible and $\Theta$-admissible $B(\Y)$-valued kernels, respectively. Let $k\in \K_\Psi(\Y)$. Then for each $\psi \in \Psi$, the map $M_\psi : \mathcal{H}(k) \rightarrow \mathcal{H}(k)$, sending $h$ to $\psi \cdot h$ is a contraction. It is a well known fact that for an $a\in \mathbb{D}$ and a contraction $T$ on a Hilbert space $\mathfrak{H}$, $(aI_{\mathfrak{H}} - T)(I_{\mathfrak{H}} -\overline{a}T)^{-1}$ is again a contraction on $\mathfrak{H}$ (see \cite{Nagy-Foias}, page 14). So if $a= \psi (w_0)$, then
$$(aI_\Y - M_\psi)(I_\Y - \overline{a}M_\psi)^{-1}: \mathcal{H}(k)\rightarrow \mathcal{H}(k)$$
is again a contraction. We shall show that $(aI_\Y - M_\psi)(I_\Y - \overline{a}M_\psi)^{-1} =M_{\theta_\psi}$. Choose any $h\in \mathcal{H}(k)$ and any $k(\cdot , w)y\in \mathcal{H}(k)$ with $w\in \Omega$ and $y\in \Y$. It is clear that $$\Big\langle(aI_\Y - M_\psi)(I_\Y - \overline{a}M_\psi)^{-1}h, k(\cdot , w)y \Big\rangle_{\mathcal{H}(k)} =(a- \psi(w))\Big\langle(I_\Y - \overline{a}M_\psi)^{-1}h, k(\cdot , w)y\Big\rangle_{\mathcal{H}(k)}.$$
Since $|a|<1$ and $M_\psi$ is a contraction, we have an absolutely convergent series representation $$(I_\Y - \overline{a}M_\psi)^{-1}=\sum_{j=0}^{\infty}\overline{a}^j M_\psi ^j $$
Each term of this series is a repeated application of $M_\psi$. Hence we have $$\Big\langle(I_\Y - \overline{a}M_\psi)^{-1}h, k(\cdot , w)y\Big\rangle_{\mathcal{H}(k)} = (1- \overline{a}\psi(w))^{-1}\Big\langle h, k(\cdot , w)y\Big\rangle_{\mathcal{H}(k)} .$$
Since $\theta_\psi (w)=\frac{a- \psi(w)}{1- \overline{a}\psi(w)}$, we get that
$$\Big\langle(aI_\Y - M_\psi)(I_\Y - \overline{a}M_\psi)^{-1}h, k(\cdot , w)y \Big\rangle_{\mathcal{H}(k)} =\theta_\psi (w) \Big\langle h, k(\cdot , w)y\Big\rangle_{\mathcal{H}(k)}.$$
This proves $\Big[(aI_\Y - M_\psi)(I_\Y - \overline{a}M_\psi)^{-1}\Big]^* = M_{\theta_\psi} ^*$ and consequently $M_{\theta_\psi}$ is a contraction on $\mathcal{H}(k)$.  Thus $k\in \K_\Theta (\Y)$. So $\K_\Psi (\Y)\subseteq\K_\Theta (\Y)$. Since $\varphi_\psi ^{-1} =\varphi_\psi $  and $\varphi_\psi ^{-1} \circ \theta_\psi = \psi$, the same argument used above proves $K_\Theta (\Y) \subseteq \K_\Psi (\Y)$. This completes the proof of the proposition.\qed

{\em From now on, we shall assume without loss of generality that our collection of holomorphic test functions has a common zero, that is, there is a point $w_0 \in \Omega$ such that $E(w_0)= 0$ in $\CB$.}

We now develope a power series representation.

Suppose $z_0\in \Omega$. Then there is a polydisc $P=P(z_0,r)$ centered at $z_0$ with multi-radius $r$ such that $\overline{P}\subset \Omega$. Let $b P$ be its distinguished boundary. Since the collection $\Psi$ consists of holomorphic functions, for each $\psi \in \Psi$ we have
\begin{align}
\psi(z)= \sum_{\alpha \in \mathbb{N}^m} \Big(\frac{1}{(2\pi i)^m} \int_{b P} \frac{\psi (\zeta)}{(\zeta - z_0)^{\alpha +1}} d\zeta\Big)(z-z_0)^\alpha
\end{align}
where $\alpha +1$ stands for $(\alpha_1 +1,\alpha_2 +1,\ldots,\alpha_m +1)$.

Now we shall introduce a collection of functions in $\CB$ that has a great importance in the power series development of $\Psi$-Schur-Agler class.
 For each $\alpha \in \mathbb{N}^m$, define a function $\delta_\alpha : \Psi \rightarrow \mathbb{C}$ by
\begin{align}\label{Definition of delta alpha}
\delta_\alpha (\psi) = \frac{1}{(2\pi i)^m} \int_{b P} \frac{\psi (\zeta)}{(\zeta - z_0)^{\alpha +1}} d\zeta .
\end{align}
We have the following result on these $\delta_\alpha$.
\begin{proposition}
	With $z_0$ and $P$ as above, $\delta_\alpha\in\CB$ for each $\alpha\in \mathbb{N}^m$.
	\end{proposition}
\textbf{Proof.}\\
Clearly $|\delta_\alpha (\psi)|\leq \frac{1}{r^\alpha}$ as $||\psi||_\Omega \leq 1$. Since this holds for any $\psi$, we have $||\delta_\alpha||\leq \frac{1}{r^\alpha}$, that is, $\delta_\alpha$ is a bounded function on $\Psi$ with norm less than or equal to $\frac{1}{r^\alpha}$. Now we shall show that $\delta_\alpha$ is continuous on $\Psi$ for each $\alpha \in \mathbb{N}^m$. Let us fix an $\alpha \in \mathbb{N}^m$ and consider a net $\{\psi_\nu :\nu \in N\}$ that converges to a $\psi$ in $\Psi$. Consider the Borel probability measure
$$d\mu =\frac{d\theta}{(2\pi)^m}=\frac{d\theta_1d\theta_2\ldots d\theta_m}{(2\pi)^m}$$ on $[0,2\pi]^m$
and a positive $\epsilon$. Let $E_\nu =\{\theta \in [0,2\pi]^m: |\psi_\nu (z_0 + re^{i\theta})-\psi (z_0 + re^{i\theta})|\geq\epsilon\}$.\\
 \textbf{Claim:} $\lim_\nu \mu (E_\nu) =0$.\\
 \textbf{Proof of the claim:}\\
 If not, then there is a positive number $\tau$ such that for every $\nu \in N$, we can find a $\eta_\nu \in N$ such that $\eta_\nu \geq \nu$ and $\mu(E_{\eta_\nu})\geq \tau$. Now $\{\psi_{\eta_\nu} :\nu \in N\}$ is a subnet of $\{\psi_\nu :\nu \in N\}$, so it converges point-wise to $\psi$ on $\Omega$. Since our collection $\Psi$ consists of holomorphic functions, we can take the closure of the uniformly bounded set $\{\psi_{\eta_\nu}:\nu \in N\}$ in the compact-open topology of the space of all holomorphic functions on $\Omega$ and get a compact subset as a result. Hence $\{\psi_{\eta_\nu}:\nu \in N\}$ has a convergent subnet $\{\psi_{\eta_\nu}:\nu \in N^\prime\}$ ($N^\prime \subset N$). Let $\{\psi_{\eta_\nu}:\nu \in N^\prime\}$ converge to some holomorphic $\phi: \Omega \rightarrow \mathbb{C}$. This convergence is uniform on each compact subset of $\Omega$. But $\psi$ is a point-wise limit of the subnet $\{\psi_{\eta_\nu}:\nu \in N^\prime\}$ as well. So we have $\phi=\psi$. Since $b P$ is a compact subset of $\Omega$, we can find a $\nu \in N^\prime$ such that
 $$\sup_{\zeta \in b P} |\psi_{\eta_\nu} (\zeta) -\psi(\zeta)|=\sup_{\theta \in [0,2\pi]^m} |\psi_{\eta_\nu} (z_0 + re^{i\theta})-\psi (z_0 + re^{i\theta})|<\epsilon.$$
 So we have $E_{\eta_\nu} =\emptyset$. But this is a contradiction, as $\mu(E_{\eta_\nu})\geq \tau>0$. Thus we get $\lim_\nu \mu (E_\nu)=0$. Hence the claim follows.\qed
 From this we conclude
 $$\lim_\nu ||\psi_\nu (z_0 + re^{i(\cdot)}) -\psi(z_0 + re^{i(\cdot)})||_{L^1 (\mu, [0,2\pi]^m)}=0$$ (see \cite{Dunford-Schwarz}
 page 124, Theorem 7). Since
 $$|\delta_\alpha(\psi_\nu) -\delta_\alpha(\psi)|\leq \frac{1}{r^\alpha} ||\psi_\nu (z_0 + re^{i(\cdot)}) -\psi(z_0 + re^{i(\cdot)})||_{L^1 (\mu, [0,2\pi]^m)},$$
 it follows that $\delta_\alpha$ is continuous on $\Psi$. So we conclude that $\delta_\alpha\in \CB$ for each $\alpha$. This completes the proof.\qed
 Now let us take a look at the maps $E(z),z\in \Omega$. The following result gives a local power series representation of $E(z)$.
 \begin{theorem}\label{Expansion of evaluation functions}
 	For each $z_0\in \Omega$, there is a polydisc $P\subset \subset \Omega$ centered at $z_0$ such that the series
 	$$\sum_{\alpha\in \mathbb{N}^m} (z-z_0)^\alpha \delta_\alpha$$
 	converges absolutely and uniformly on every compact subset of $P$ to $E(z)$.
 	\end{theorem}
 \textbf{Proof.}\\
 Take $P$ and $r$ as above and fix a $z\in P$. Then $$|z-z_0|^\alpha ||\delta_\alpha||\leq \frac{|z-z_0|^\alpha}{r^\alpha}$$
  for every $\alpha \in \mathbb{N}^m$. Since $z$ is in the interior of $P$, we can find a $q\in (0,1)$ such that $\frac{|z-z_0|^\alpha}{r^\alpha}<q^{|\alpha|}$ for each $\alpha$. So the series $\sum_{\alpha \in \mathbb{N}^m} |z-z_0|^\alpha ||\delta_\alpha||$ converges. It also gives that $\sum_{\alpha\in \mathbb{N}^m} (z-z_0)^\alpha \delta_\alpha$ uniformly on the closed polydisc $\overline{P}(z_0, |z-z_0|)$. This settles the requirements of absolute convergence and uniform convergence on compact subsets. Thus the series in question does define an element of $\CB$.

 Now we have for any $\psi\in\Psi$
 \begin{align*}
 \psi(z)= \sum_{\alpha \in \mathbb{N}^m} \Big(\frac{1}{(2\pi i)^m} \int_{b P} \frac{\psi (\zeta)}{(\zeta - z_0)^{\alpha +1}} d\zeta\Big)(z-z_0)^\alpha = \sum_{\alpha \in \mathbb{N}^m} \delta_\alpha (\psi) (z-z_0)^\alpha.
 \end{align*}
 So if $F\subset \mathbb{N}^m$ is a finite set and if $z^\prime \in\overline{P}(z_0, |z-z_0|)$, then
 \begin{align*}
 |E(z^\prime)(\psi)- \sum_{\alpha \in F} \delta_\alpha (\psi) (z^\prime -z_0)^\alpha|\leq \sum_{\alpha \notin F}\frac{|z-z_0|^\alpha}{r^\alpha}\leq \sum_{\alpha \notin F} q^{|\alpha|}
 \end{align*}
with $q$ as above. Clearly the convergence of the rightmost expression is not affected by our choice of $z^\prime$ from $\overline{P}(z_0, |z-z_0|)$. Hence, $E(z)=\sum_{\alpha\in \mathbb{N}^m} (z-z_0)^\alpha \delta_\alpha$ for all $z\in P$ and the convergence is uniform on every compact subset of $\Omega$. This completes the proof.\qed
Let $\X$ be any Hilbert space and $\rho:\CB \rightarrow B(\X)$ a unital $*$-representation. Then $||\rho||=1$ and consequently we are allowed to write
$$\rho(E(z))=\sum_{\alpha\in \mathbb{N}^m} (z-z_0)^\alpha \rho(\delta_\alpha)$$
for all $z\in P$. The consequences of Theorem \ref{Expansion of evaluation functions} do not cease here. Absolute convergence allows us to use the same procedure that is used to prove Mertens' theorem (see Theorem $3.4$ in \cite{Baby-Rudin}) and conclude that for any integer $l\geq 1$ and a contraction $A:\X\rightarrow\X$ the map $z\mapsto\rho(E(z)) A \rho(E(z))^l $ has a power series representation on $P$ such that the series converges absolutely and uniformly on every compact subset of $P$.

The power series representation also implies Bochner integrability in some smaller sets. Consider a Banach space  $\B$ and a map $f:P(z_0,r)\rightarrow \B$ with a power series representation
$$f(z)=\sum_{\alpha \in \mathbb{N}^m} (z-z_0)^\alpha a_\alpha$$
where $a_\alpha \in \B$. The convergence of this series is absolute and uniform on every compact subset of $P(z_0,r)$. We take a smaller polydisc $P^\prime = P^\prime (z_0,r^\prime)\subset \subset P(z_0,r)$. Then $f(b P^\prime)$ is separable, because the collection $\big\{\sum_{\alpha\in F} (\zeta -z_0)^\alpha a_\alpha : F (\subset\mathbb{N}^m)\,\, \text{is finite and}\,\, \zeta\in \mathbb{Q}^m + i\mathbb{Q}^m\big\}$ is dense in $f(b P^\prime)$. Also $f$ is continuous. Hence $f$ is strongly Borel measurable (see \cite{Cohn}). Now the absolute and uniform convergence of the power series on $b P^\prime$ gives us that for any finite Borel measure $\mu$ on $b P^\prime$
$$\int_{b P^\prime} f(\zeta) d\mu (\zeta)= \sum_{\alpha \in \mathbb{N}^n} \Big(\int_{b P^\prime} (\zeta -z_0)^\alpha d\mu(\zeta) \Big) a_\alpha.$$
From this discussion, it follows that for a function $f$ with a power series representation, we have a Cauchy like formula, that is,
$$f(z)=\frac{1}{(2\pi i)^m}\int_{b P^\prime} \frac{f(\zeta)}{(\zeta -z)} d\zeta$$
for all $z\in P^\prime$. Obviously a power series representation on $P^\prime$ can be deduced from this integral formula.

Recall the polydisc $P(z_0,r)$ that was used to develop the power series representation of $E(z)$. Since $P(z_0,r)\subset \subset \Omega$, a larger multi-radius $r$ could have been considered without any difficulty in dealing with the power series. So we may assume that the map $\mathfrak{E}:\overline{P}(z_0,r) \rightarrow \mathbb{R}$ sending $z$ to $||E(z)||$ is continuous and hence attains a maximum at some point of $\overline{P}(z_0,r)$. Since the maximum is attained and $||E(z)||<1$ for each $z\in \Omega$, we can find a $q_0 \in (0,1)$ such that $\sup_{\overline{P}(z_0,r)}\mathfrak{E} <q_0$. This $q_0$ will play an important role in our next discussion.

Let us take $P(z_0, r), E(z), \X, \rho$ and $A$ as above. For each $l\in\mathbb{N}$ define a function $g_l :P(z_0, r) \rightarrow B(\X)$ by
$$g_l (z)=\rho(E(z)) A \rho(E(z))^l .$$
Clearly $||g_l (z)||< q_0 ^{l+1}$ whenever $z\in P(z_0, r)$. So the series $h(z)= \sum_{l\geq 1} g_l (z)$ converges absolutely and uniformly on every compact subset of $P(z_0, r)$. Now each $g_l$ has a power series representation in $P(z_0, r)$ and hence Bochner integrable on the distinguished boundary of some smaller polydisc $P^\prime$ centered at $z_0$. Consequently, $h$ is also Bochner integrable on the same set. Putting all this together, we get that
$$h(z)=\frac{1}{(2\pi i)^m}\int_{b P^\prime} \frac{h(\zeta)}{(\zeta -z)} d\zeta$$
for all $z\in P^\prime$. Using classical analysis, we conclude that $h$ has a power series representation in some neighborhood of $z_0$. Since we can get a series expansion in $\rho(E(z))$ from the realization formula for $S\in \SA_\Psi (\U,\Y)$ and the series takes the same form as $h$ above, we obtain the following theorem.

\begin{theorem}\label{Power series for Schur-Agler class}
	Let $\Omega$ be a domain in $\mathbb{C}^m$ and let $\Psi$ be a class of holomorphic test functions on $\Omega$. Then for each $z_0 \in \Omega$ and $S\in \SA_\Psi(\U,\Y)$, there is a polydisc $P(z_0, \varepsilon)\subset \subset\Omega $ centered at $z_0$ with multi-radius $\varepsilon$ such that
	\begin{align*}
	S(z)= \sum_{\alpha \in \mathbb{N}^m} (z-z_0)^\alpha a_\alpha\,\,\text{for all}\,\, z\in P(z_0, \varepsilon),
	\end{align*}
	where $a_\alpha \in B(\U,\Y)$. Moreover, the series converges absolutely and uniformly on each compact subset of $P(z_0, \varepsilon)$.
\end{theorem}
\qed

\section{An Auxiliary Function}\label{Affiliated}

We start with a given collection $\Psi$ of test functions  that vanish at the point $w_0$. So $E(w_0)=0$ in $\CB$. Let us prove the following lemma. \\

\textbf{Lemma.} Suppose that $z_i \in \Omega$ and  $B_i \in B(\U, \Y)$, $1 \leq i \leq n$. The interpolation problem $z_i \mapsto B_i$ is solvable by a function in $\SA_\psi (\U, \Y)$ if and only if there is a completely positive kernel $\Gamma: \{z_1,z_2,\ldots, z_n\} \times \{z_1,z_2,\ldots, z_n\}\rightarrow B(\mathcal{C}_b(\Psi), B(\Y))$ such that
	\begin{align}\label{DataDecomp}
	I_\Y - B_i B_j ^*  & = \Gamma(z_i, z_j)(1- E(z_i) E(z_j)^*) \ \ \ \text{for all} \ \ i,\, j =1,2,\ldots , n.
	\end{align}
The proof of this lemma follows from Theorem \ref{RealTheo}.\qed

Now let us construct a function $G$. For a solvable problem $z_i \mapsto B_i, 1\leq i\leq n$, (\ref{DataDecomp}) gives us a completely positive kernel $$\Gamma: \{z_1,z_2,\ldots, z_n\} \times \{z_1,z_2,\ldots, z_n\}\rightarrow B(\mathcal{C}_b(\Psi), B(\Y))$$ such that
\begin{align*}
I_\Y - B_i B_j ^*  & = \Gamma(z_i, z_j)(1- E(z_i) E(z_j)^*) \ \ \ \text{for all} \ \ i,\, j =1,2,\ldots , n.
\end{align*}
 Using Kolmogorov decomposition (\ref{Kolmogorv Decomposition}) for $$\Gamma: \{z_1,z_2,\ldots, z_n\} \times \{z_1,z_2,\ldots, z_n\} \rightarrow B(\mathcal{C}_b(\Psi), B(\Y))$$  we have that there exists a Hilbert space $\X_1$, a unital $*$-representation $\mu: \mathcal{C}_b(\Psi) \rightarrow B(\X_1)$ and a function $h: \{z_1,z_2,\ldots, z_n\} \rightarrow B(\X_1, \Y)$ such that
\begin{align*}
\Gamma(z_i, z_j)(\delta) = h(z_i) \mu(\delta) h(z_j)^* \,\,\text{for all} \ \ i, j \in \{1,2,\ldots, n\} \ \ \text{and} \ \ \delta \in \mathcal{C}_b(\Psi).
\end{align*}
This gives us
\begin{align*}
I_\Y - B_i B_j ^*  & = h(z_i)h(z_j)^* - h(z_i)\mu (E(z_i) E(z_j)^*)h(z_j)^*.
\end{align*}
So for any $y,w \in \Y$, we have
\begin{align*}
\langle y,w\rangle - \langle B_i ^* y, B_j ^* w\rangle = \langle h(z_i)^*y, h(z_j)^* w\rangle -  \langle \mu (E(z_i))^* h(z_i)^*y, \mu (E(z_j))^*  h(z_j)^* w\rangle .
\end{align*}
Let $\mathscr{L}_1=\overline{span}\{\mu (\delta)h(z_i)^* y: 1\leq i\leq n, \delta \in \mathcal{C}_b(\Psi),y\in \Y\}$. So the last equality can be rewritten as
\begin{align*}
\Bigg\langle \begin{pmatrix} \mu(E(z_i))^* h(z_i)^* y  \\ y \end{pmatrix},  \begin{pmatrix}  \mu(E(z_j))^* h(z_j)^* w\\ w  \end{pmatrix} \Bigg\rangle_{ \mathscr{L}_1\oplus \Y} = \Bigg\langle \begin{pmatrix}  h(z_i)^* y\\ B_i^* y \end{pmatrix}, \begin{pmatrix} h(z_j)^* w\\ B_j^* w  \end{pmatrix} \Bigg\rangle_{ \mathscr{L}_1\oplus \U}.
\end{align*}
Now, let
$$\N_2 = \overline{span} \Bigg\{   \begin{pmatrix} \mu(E(z_i))^* h(z_i)^* y \\  y\end{pmatrix} : y \in \Y, i = 1,2, \dots ,n \Bigg\}$$
and
$$\N_1 = \overline{span} \Bigg\{   \begin{pmatrix}  h(z_i)^* y \\ B_i^* y \end{pmatrix} : y \in \Y, i = 1,2, \dots ,n \Bigg\}.$$
Then there is a unitary $V:\mathscr{N}_2\rightarrow \mathscr{N}_1$ that sends

\begin{align}\label{V}
\begin{pmatrix} \mu(E(z_i))^* h(z_i)^* y  \\ y \end{pmatrix}\,\, \text{to}\,\, \begin{pmatrix}  h(z_i)^* y\\ B_i^* y \end{pmatrix}\,\, \text{for all}\,\, i\,\, \text{and}\,\, y.
\end{align}

Let $\M_1 = ( \mathscr{L}_1 \oplus \U) \ominus \N_1$ and $\M_2 = ( \mathscr{L}_1\oplus \Y) \ominus \N_2$. Since, $\N_1$ and $\N_2$ are unitarily equivalent, the linear operator $\Q : \N_2 \oplus \M_2 \oplus \M_1 \rightarrow \N_1 \oplus \M_1 \oplus \M_2$ sending  \ $n_2 \oplus m_2 \oplus m_1$ to $Vn_2 \oplus m_1 \oplus m_2$ is a well defined unitary operator. Since $\N_2 \oplus \M_2 \oplus \M_1 \simeq  \mathscr{L}_1 \oplus \Y \oplus \M_1$ and $\N_1 \oplus \M_1 \oplus \M_2 \simeq \mathscr{L}_1 \oplus \U \oplus  \M_2$, we can write $\Q$ as
\begin{align}\label{Defn-Q}
\Q  =
\bordermatrix{ & \mathscr{L}_1 & \M_1 \oplus \Y  \cr
	\mathscr{L}_1 & \Q_{11} & \Q_{12} \cr
	\M_2 \oplus \U & \Q_{21} & \Q_{22} }. \qquad
\end{align}
  We consider the function
\begin{align}\label{defn_G}
G(z) := \Q^*_{22} + \Q^*_{12} (I_{\mathscr{L}_1} - \mu(E(z)) \ \Q_{11}^* )^{-1} \  \mu(E(z)) \ \Q_{21}^*, \ (z\in \Omega).
\end{align}
Then for each $z \in \Omega$, we have $G(z) \in B(\M_2 \oplus \U, \M_1 \oplus \Y)$. Moreover, by the result of Section \ref{Char-schur-agler}, we have $G \in \SA_\psi(\M_2 \oplus \U, \M_1 \oplus \Y)$. Since, $\Q$ is an extension of $V$, we have that $$\Q(n_2 \oplus \mathbf{0}_{\M_2} \oplus \mathbf{0}_{\M_1}) = V{n_2} \oplus \mathbf{0}_{\M_1} \oplus \mathbf{0}_{\M_2} \ \ \text{for all} \ \ n_2 \in \N_2.$$\\
\noindent
Hence, we can write
\begin{align*}
\Q \begin{pmatrix}  \mu(E(z_i))^* h(z_i)^* y\\  y  \end{pmatrix}  & =  \begin{pmatrix}  h(z_i)^* y\\  B_i^* y      \end{pmatrix}
\end{align*}
for all $ i = 1,2, \dots ,n $ and $ y \in \Y$. In other words,
\begin{align*}
\bordermatrix{ & \mathscr{L}_1 & \M_1 \oplus \Y  \cr
	\qquad \mathscr{L}_1 & \Q_{11} & \Q_{12} \cr
	\M_2 \oplus \U & \Q_{21} & \Q_{22} } \begin{pmatrix} \mu(E(z_i))^* h(z_i)^* y \\ \mathbf{0}_{\M_1} \oplus y \end{pmatrix} & = \begin{pmatrix} h(z_i)^* y \\ \mathbf{0}_{\M_2} \oplus B_i^* y \end{pmatrix},
\end{align*}
\begin{align*}
\text{or,} \ \  \Q_{11} (\mu (E(z_i))^* h(z_i)^* y) + \Q_{12}( \mathbf{0}_{\M_1} \oplus y)  & = h(z_i)^* y \\
\text{and} \ \ \Q_{21}(\mu (E(z_i))^* h(z_i)^* y) + \Q_{22}( \mathbf{0}_{\M_1} \oplus y)  & = \mathbf{0}_{\M_2} \oplus B_i^* y,
\end{align*}
\begin{align*}
\text{or,} \ \  (I_{\mathscr{L}_1} - \Q_{11} \mu(E(z_i))^*)^{-1} \mathcal{Q}_{12}( \mathbf{0}_{\M_1} \oplus y)  & = h(z_i)^* y \\
\text{and} \ \ \Q_{22}( \mathbf{0}_{\M_1} \oplus y) + \Q_{21}(\mu (E(z_i))^* h(z_i)^* y) & = \mathbf{0}_{\M_2} \oplus B_i^* y.
\end{align*}
Combining the last two equations, we obtain
\begin{align*}
\Q_{22}(\mathbf{0}_{\M_1} \oplus y) + \Q_{21}\mu(E(z_i))^* & (I_{\mathscr{L}_1} - \Q_{11} \mu (E(z_i))^*)^{-1} \Q_{12}(\mathbf{0}_{\M_1} \oplus y)=\mathbf{0}_{\M_2} \oplus B_i^* y,
\end{align*}
and hence from the definition of $G$ (see \eqref{defn_G}) we get
\begin{align}\label{G}
G(z_i)^* (\mathbf{0}_{\M_1} \oplus y)= \mathbf{0}_{\M_2} \oplus B_i^* y.
\end{align}
Let us write $G(z)^*$ as
\begin{align*}
G(z)^* = \bordermatrix{ & \M_1 &  \Y  \cr
	\M_2 & G_{11}(z_i)^* & G_{21}(z_i)^* \cr
	\U & G_{12}(z_i)^* &G_{22}(z_i)^* }.
\end{align*}
So we have
\begin{align*}
\bordermatrix{ & \M_1 &  \Y  \cr
	\M_2 & G_{11}(z_i)^* & G_{21}(z_i)^* \cr
	\U & G_{12}(z_i)^* &G_{22}(z_i)^* } \begin{pmatrix} \mathbf{0}_{\M_1}  \\ y \end{pmatrix} =  \begin{pmatrix} \mathbf{0}_{\M_2}  \\ B_i^* y \end{pmatrix}
\end{align*}
which gives
\begin{align*}
G_{21}(z_i) =  \mathbf{0} \ \ & \text{and} \ \ G_{22}(z_i) = B_i \ \ \text{for} \ \ i = 1,2, \dots ,n.
\end{align*}
This completes the construction of $G$.

Note that for any $m_1 \in \M_1$, we have $\Q(\mathbf{0}_{\N_2}\oplus \mathbf{0}_{\M_2}\oplus m_1)=(\mathbf{0}_{\N_1}\oplus m_1 \oplus \mathbf{0}_{\M_2})$. Since $\N_2 \oplus \M_2=\LL_1\oplus \Y$ and $\N_1 \oplus \M_1=\LL_1\oplus \U$, we can write $\mathbf{0}_{\N_2}\oplus \mathbf{0}_{\M_2} = \mathbf{0}_{\LL_1} \oplus \mathbf{0}_{\Y}$ and $\mathbf{0}_{\N_1}\oplus m_1= l\oplus u$ for some $l\in \LL_1$ and $u\in \U$. So using (\ref{Defn-Q}) we get that $\Q_{22} (m_1 \oplus \mathbf{0}_\Y)=\mathbf{0}_{\M_2}\oplus u$. Now we have that for any $z\in \Omega$,
$$G(z)^* := \Q_{22} + \Q_{21} (I_{\mathscr{L}_1} - \mu(E(z))^* \ \Q_{11} )^{-1} \  \mu(E(z))^* \ \Q_{12};$$
and by virtue of Proposition \ref{Modified Test Functions}, we can assume without loss of generality that there is a point $w_0\in\Omega$ such that $E(w_0)=0$. So we have $G(w_0)^* (m_1 \oplus \mathbf{0}_\Y)=\mathbf{0}_{\M_2}\oplus u$. Since $m_1 \in \M_1$ is arbitrary, this gives $G_{11}(w_0) \equiv \mathbf{0}$ as an operator on $\M_2$. Hence the function $z\mapsto ||G_{11} (z)||$ from $\Omega$ to $\mathbb{R}$ is non-constant. Being a component of a $\Psi$-Schur-Agler class function $G$, the function $G_{11}$ has a power series expansion by Theorem \ref{Power series for Schur-Agler class}. So we can apply maximum modulus theorem for Banach space valued holomorphic functions and deduce that if $z\in \Omega$, then $||G(z)||<1$. A proof of the maximum modulus theorem for Banach space valued holomorphic functions of one variable can be found in (\cite{Taylor-Lay}, page 269, Theorem $1.5$), and this proof can easily be carried out in our case as well.

For a given set of data for a solvable interpolation problem, $z_i\in \Omega$ and $B_i\in \B(\U,\Y)$ $1\leq i\leq n$, we shall keep this $G$ fixed.

\section{The Main Result}\label{Main Result}

We start this section by noting that given a data set $\{z_1, z_2, \ldots , z_n\}$ and $\{B_1, B_2, \ldots , B_n\}$ where the $B_i$ are in $B(\U, \Y)$, if there are two Hilbert spaces $\M_1$ and $\M_2$ and a function
$$G = \left(
        \begin{array}{cc}
          G_{11} & G_{12} \\
          G_{21} & G_{22} \\
        \end{array}
      \right)
 \text{ in } \SA_\Psi(\M_2\oplus \U, \M_1 \oplus \Y)$$
 satisfying $G_{22}(z_i) = B_i$, $G_{21}(z_i) = 0$ for all $i=1,2, \ldots , n$ and $\| G_{11}(z)\| < 1 $ at all points of $\Omega$, then for any function $\mathbf t$ in  $\SA_\Psi(\M_1 , \M_2 )$, the function
 $$f(z) = \left(G_{22} + G_{21}(I_{\M_2} - \mathbf t G_{11})^{-1} \mathbf t G_{12}\right)(z)$$
 is an interpolant for the data in $\SA_\Psi( \U,  \Y)$.

Consider a solvable interpolation problem $z_i \mapsto B_i$, where $z_i \in \Omega$ and $B_i \in B(\U, \Y)$, $1 \leq i \leq n$. If $f\in  \SA_\Psi(\U, \Y)$ is a solution to this interpolation problem then Theorem \ref{RealTheo} gives us that we can find a completely positive kernel $\Delta: \Omega \times \Omega\rightarrow B(\mathcal{C}_b(\Psi), B(\Y))$ such that
\begin{align}\label{SolDecomp}
I_{\Y} - f(z) f(w)^* = \Delta(z,w) (1- E(z) E(w) ^*) \ \ \text{for all} \ z, w \in \Omega.
\end{align}
The restriction of  $\Delta$ to $\{z_1, z_2, \ldots , z_n\} \times \{z_1, z_2, \ldots , z_n\}$ may not agree with $\Gamma$ in (\ref{DataDecomp}) in general. When they do, we call $f$ an \textit{affiliated solution}. To be more precise, we give a proper definition.

\textbf{Definition}: Let $z_1, z_2, \ldots, z_n\in \Omega$ and $B_1, B_2, \ldots, B_n \in B(\U, \Y)$ be a solvable data. Let $f \in \SA_\psi(\U, \Y)$ be a solution. Let $\Gamma$ and $\Delta$ be as in (\ref{DataDecomp}) and (\ref{SolDecomp}), respectively. Then $f$ is said to be affiliated with $\Gamma$ if $\Gamma(z_i, z_j)= \Delta(z_i, z_j)$ for all $i,j=1,2,\ldots , n$.

Why does one need the concept of affiliation? Because, given a solvable interpolation problem, the kernel $\Gamma$ obtained in \eqref{DataDecomp} may not be unique. An example can be found on page $185$ of \cite{A-M}. The notion of affiliation first appeared in \cite{B-T} where the authors solved the Nevanlinna problem for the bidisc assuming this notion. The following theorem is a generalization.

\begin{theorem} \label{parametrization}
	Let $\Omega$ be a domain in $\mathbb{C}^m$ and let $\Psi$ be a class of holomorphic test functions. Suppose that $f \in \SA_\Psi(\U, \Y)$ is a solution of this interpolation problem $z_i \mapsto B_i$ and $f$ is affiliated with a completely positive kernel
	\begin{align*}
	\Gamma: \{z_1,z_2,\ldots, z_n\} \times  \{z_1,z_2,\ldots, z_n\}\rightarrow B(\mathcal{C}_b(\Psi), B(\Y)).
	\end{align*}
	Let $\M_1$, $ \M_2$ and $G$ be as in Section \ref{Affiliated}. Writing $G$ as
	\[
	G(z) =
	\bordermatrix{ & \mathscr{M}_2 & \U \cr
		\mathscr{M}_1 & G_{11}(z) & G_{12}(z) \cr
		\Y & G_{21}(z) & G_{22}(z)}, \qquad
	\]
	we have $$f(z) = \Big(G_{22} + G_{21}(I_{\M_2} -\mathbf{t}\  G_{11})^{-1} \mathbf{t}\  G_{12}\Big)(z)$$ for some $\mathbf{t} \in \SA_\Psi(\M_1, \M_2)$ and for all $ z \in \Omega.$
\end{theorem}

\vspace*{5mm}

{\bf Note}: Before we embark on the proof, we want to remark that

\begin{enumerate}
\item without loss of generality, we can assume that all test functions vanish at a certain point $w_0$, i.e., $E(w_0) = 0$ because of Proposition \ref{Modified Test Functions},
    \item the inverse of $I_{\M_2} -\mathbf{t}\  G_{11}$ exists because of the concluding remarks of Section \ref{Affiliated}.
\end{enumerate}

\vspace*{5mm}

 \textbf{Proof of Theorem \ref{parametrization}.}\\
 Since $f$ is a solution and $f$ is affiliated with $\Gamma$, we can find a completely positive kernel $\Delta: \Omega \times \Omega\rightarrow B(\mathcal{C}_b(\Psi), B(\Y))$ such that
 \begin{align}\label{F-Agler-Decomp}
 I_{\Y} - f(z) f(w)^* = \Delta(z,w) (1- E(z) E(w) ^*) \ \ \text{for all} \ z, w \in \Omega
 \end{align}
and $\Gamma(z_i, z_j)= \Delta(z_i, z_j)$ for all $i,j=1,2,\ldots , n$. Now, we know from (\ref{Kolmogorv Decomposition}) that there is a Hilbert space $\mathscr{X}$,  $*$-representation $\rho: \mathcal{C}_b(\Psi) \rightarrow B(\X)$ and a function $g: \Omega \rightarrow B(\X, \Y)$ such that
\begin{align} \label{rep}
	\Delta(z, w)(a) = g(z) \rho(a) g(w)^* \,\,\text{for all} \ \ z, w \in \Omega \ \ \text{and} \ \ a \in \mathcal{C}_b(\Psi).
	\end{align}
 From these equations we can construct a unitary $\tilde{W}:\mathscr{X}\oplus \Y \rightarrow \mathscr{X}\oplus \U$ (as we did in section \ref{Affiliated}) such that writing $\tilde{W}$ as
	\[
\tilde{W} =
\bordermatrix{ & \X & \Y \cr
	\X & \tilde{A} & \tilde{B} \cr
	\U & \tilde{C} & \tilde{D}}, \qquad
\]
one has
\begin{align}\label{FandWtilde}
f(z)^* = \tilde{D}+ \tilde{C} (I_{\X}- \rho(E(z))^* \tilde{A}) ^{-1} \rho(E(z))^* \tilde{B} \ \ \text{for all} \ z \in \Omega
\end{align}
and $\tilde{W}$ takes $\begin{pmatrix} \rho (E(z))^* g(z)^* y \\  y \end{pmatrix} \text{ to } \begin{pmatrix} g(z)^* y \\ f(z)^* y \end{pmatrix}.$\\

Let $$\LL= \overline{span}\{\rho(\delta) g(z_i)^* y: y\in \Y, \delta \in \CB, 1\leq i\leq n \}.$$ This is a closed subspace of $\X$ and it is reducing for $\rho(E(z)),\, \text{for all}\, z\in \Omega$. Recall the subspace $\LL_1$ of $\X_1$ from Section \ref{Affiliated} which was defined by $$\mathscr{L}_1=\overline{span}\{\mu (\delta)h(z_i)^* y: 1\leq i\leq n, \delta \in \mathcal{C}_b(\Psi),y\in \Y\}.$$ Now $\Gamma(z_i, z_j)= \Delta(z_i, z_j)$ for all $i,j=1,2,\ldots , n,$ gives us $$ g(z_i)\rho (\delta)g(z_j)^* =h(z_i)\mu (\delta)h(z_j)^*,\, \text{for all}\, 1\leq i,j\leq n,\, \delta \in \CB.$$ It is easy to see that the map $\tilde{S}:\LL\rightarrow \LL_1$ sending $\rho(\delta) g(z_i)^* y$ to $\mu (\delta)h(z_i)^* y$ is a unitary.\\
Let $\HH = \X \ominus \LL$ and $S =I_{\HH} \oplus \tilde{S}$. Then $S:\mathscr{X} \rightarrow \HL$ is a unitary.\\
We define $\lambda: \CB \rightarrow B(\HL)$ by $$\lambda (\delta) = S\rho(\delta)S^*,\,\delta\in \CB,$$ and $l:\Omega \rightarrow B(\HL \oplus \Y)$ by $$l(z)= g(z)S^*.$$
Clearly, $\lambda$ is a unital $*$-representation. We have
\begin{align}\label{Defn-lambda}
l(z_i)^* y = h(z_i)^* y\,\, \text{and}\,\, \lambda(\delta)l(z_i)^* y = \mu(\delta)h(z_i)^* y\,\, \text{for all}\,\, 1\leq i,j\leq n,\, \delta \in \CB.
\end{align}
So for any $\delta, \delta^\prime\in \CB$, $y\in\Y$ and $i=1,2,\ldots,n$, we get  $\lambda(\delta)\big(\mu(\delta^\prime)h(z_i)^*y\big)=\mu(\delta)\big(\mu(\delta^\prime)h(z_i)^*y\big)$, that is,
\begin{align}\label{equality of lambda and mu}
\lambda (\delta)|_{\LL_1} = \mu(\delta)|_{\LL_1},\,\, \text{for all}\,\, \delta \in \CB.
\end{align}
Since $\LL_1$ is reducing for $\lambda(E(z))$ for all $z\in \Omega$, we get using (\ref{rep}) and the definitions of $\lambda$ and $l$ that
\begin{align*}
I_\Y - f(z) f(w)^* = l(z) \lambda(1- E(z) E(w)^*) l(w)^*\,\, \text{for all}\,\, z,w\in \Omega.
\end{align*}
Write $W= (S\oplus I_\U ) \tilde{W}(S\oplus I_\Y)^*$. This is a unitary from $\HL \oplus \Y$ to $\HL\oplus\U$ that takes $$\begin{pmatrix} \lambda (E(z))^* l(z)^* y \\  y \end{pmatrix} \text{ to } \begin{pmatrix} l(z)^* y \\ f(z)^* y \end{pmatrix} \text{for all}\,\, y\in \Y\,\text{and}\,\, z\in \Omega.$$
So writing $W$ as
\[
W =
\bordermatrix{ & \HL & \Y \cr
	\HL & A & B \cr
	\U & C & D}, \qquad
\]
one has
\begin{align}\label{f-colli}
f(z)^* = D + C \lamZ^* (I_{\HL}- A\lamZ ^*)^{-1} B ,\,\, \text{for all}\,\, z\in \Omega.
\end{align}
In particular, $W$ takes
\begin{align}\label{WV}
\begin{pmatrix} \mu (E(z_i))^* h(z_i)^* y \\  y \end{pmatrix} \text{ to } \begin{pmatrix} h(z_i)^* y \\ B_i ^* y \end{pmatrix} \text{for all}\,\, y\in \Y\,\text{and}\,\, 1\leq i\leq n
\end{align}
because of (\ref{Defn-lambda}).

Now recall that $\LL_1 \oplus \Y =\N_2 \oplus \M_2$ and $\LL_1 \oplus \U = \N_1 \oplus \M_1$ from Section \ref{Affiliated}. So $W$ maps $\HH \oplus\M_2 \oplus\N_2$ onto $\HH\oplus\M_1 \oplus \N_1$ and maps $\N_2$ onto $\N_1$ unitarily. So from (\ref{V}) and (\ref{WV}) we obtain $W|_{\N_2} = V$. Hence we are allowed to write

 \begin{align}\label{WandX}
W = \bordermatrix{ & \HH \oplus\M_2 &  \N_2  \cr
	\qquad \HH \oplus \M_1 & Z  & \mathbf{0} \cr
\qquad	\N_1 & \mathbf{0} & V }
\end{align}
where $Z:  \HH \oplus\M_2 \rightarrow  \HH \oplus \M_1$ is a unitary. Now we write $Z$ as
\begin{align}
Z = \bordermatrix{ & \HH &  \M_2  \cr
	\quad  \HH & Z_{11}  & Z_{12} \cr
	\quad \M_1 & Z_{21} & Z_{22} }
\end{align}
and take
 \begin{align}\label{function-t}
\mathbf{t}(z)^*  = Z_{22} +  Z_{21} \lamZ^* (I_{\HH} - Z_{11}\lamZ^* )^{-1} Z_{12}\,\, \text{for all}\,\, z\in \Omega.
\end{align}
Clearly $\mathbf{t}\in \SA_\Psi(\M_1, \M_2)$.

 Let us fix a $z\in \Omega$ and a $y\in \Y$, and put $u=f(z)^* y$.
Let $k_z = (I_{\HL}- A\lamZ ^*)^{-1} By$. It is an element of $\HL$. A little computation gives
\begin{align*}
A \lamZ^* k_z + By=k_z\,\, \text{and}\,\, C\lamZ^* k_z + Dy = u.
\end{align*}
This can be rewritten as
\begin{align}\label{W-Sol-send}
W \begin{pmatrix} \lamZ^* k_z \\ y \end{pmatrix} \  =  \begin{pmatrix} k_z \\ u \end{pmatrix}.
\end{align}

For any Hilbert space $X$ and a closed subspace $M$ of $X$, let us denote the orthogonal projection of $X$ onto $M$ by $P_{X\rightarrow M}.$

Let $r_z = P_{\HL \rightarrow \HH} k_z$ . Since $\LL_1$ is reducing for $\lamZ\,\, \text{for all}\,\, z\in \Omega$ (see (\ref{Defn-lambda})), we have the following
$$
\lamZ ^* r_z = P_{\HL \rightarrow \HH} (\lamZ ^* k_z)
=  P_{\HL \oplus \Y  \rightarrow \HH} (\lamZ ^* k_z\oplus y)\,\, \text{and}$$
 $$r_z = P_{\HL \oplus \U  \rightarrow \HH} ( k_z\oplus u).$$\\
 Since $$ \lamZ ^* k_z \oplus y \in \HL \oplus\Y= \HH \oplus \N_2\oplus\M_2\,\,\, \text{and}$$
 $$ k_z\oplus u\in \HL \oplus \U = \HH \oplus \N_1\oplus \M_1,$$
there exist $n^\prime _i \in \N_i$ and $m^\prime _i \in \M_i,\,\, i=1,2$, such that \\
\begin{align}\label{r_z +n_1 + n_2 = .. etc}
  \lamZ ^* k_z \oplus y =\lamZ ^* r_z \oplus n^\prime _2  \oplus m^\prime _2\,\,\, \text{and}\,\,
k_z \oplus u= r_z \oplus n^\prime _1\oplus m^\prime _1.
\end{align}
 So (\ref{W-Sol-send}) gives us
 \begin{align*}
 W(\lamZ ^* r_z \oplus n^\prime _2  \oplus m^\prime _2)=(r_z \oplus n^\prime _1\oplus m^\prime _1)
 \end{align*}
 and using (\ref{WandX}) we get
 \begin{align*}
 \bordermatrix{ & \HH \oplus\M_2 &  \N_2  \cr
 	\qquad \HH \oplus \M_1 & Z  & \mathbf{0} \cr
 	\qquad	\N_1 & \mathbf{0} & V }  \begin{pmatrix}\lamZ ^* r_z \oplus m^\prime _2 \\ n^\prime _2 \end{pmatrix} = \begin{pmatrix} r_z \oplus m^\prime _1 \\ n^\prime _1 \end{pmatrix}.
 \end{align*}
So
\begin{align}\label{Vn_2 =n_1}
Z(\lamZ ^* r_z \oplus m_2)=  r_z \oplus m^\prime _1\,\, \text{and}\,\,
Vn^\prime _2 =n^\prime _1.
\end{align}
Using the decomposed form of $Z$ with respect to its domain and range we get \\
\begin{align*}
 \bordermatrix{ & \HH &  \M_2  \cr
	\quad  \HH & Z_{11}  & Z_{12} \cr
	\quad \M_1 & Z_{21} & Z_{22} } \begin{pmatrix}\lamZ ^* r_z \\ m^\prime _2 \end{pmatrix} = \begin{pmatrix} r_z  \\ m^\prime _1 \end{pmatrix}.
\end{align*}
This gives us two equations from which we eliminate $r_z$. Recalling the definition of $\mathbf{t}(z)^*$ (\ref{function-t}) enables us to obtain
\begin{align}\label{t(z)m_2 =m_1}
\mathbf{t}(z)^* m^\prime _2=m^\prime _1.
\end{align}
Now let $q_z = P_{\HL \rightarrow \LL_1} (k_z).$ So from (\ref{r_z +n_1 + n_2 = .. etc}) we get
\begin{align*}\label{n_1 + m_1}
n^\prime _1 \oplus m^\prime _1 &= P_{\HH \oplus \N_1 \oplus \M_1 \rightarrow \N_1 \oplus \M_1} (r_z \oplus n_1\oplus m_1)\\ &= P_{\HL \oplus \U \rightarrow \LL_1 \oplus \U} (k_z\oplus u) = q_z \oplus u
\end{align*}
and
\begin{align*}
n^\prime _2\oplus m^\prime _2 &= P_{\HH \oplus\N_2 \oplus \M_2 \rightarrow \N_2 \oplus \M_2} (\lamZ ^* r_z \oplus n^\prime _2  \oplus m^\prime _2)\\
&= P_{\HL \oplus \Y\rightarrow \LL_1 \oplus \Y} (\lamZ^* k_z \oplus y) = \lamZ^* q_z \oplus y.
\end{align*}
Recall the $\Q$ that was defined in (\ref{Defn-Q}) in Section \ref{Affiliated}. It is a unitary from $\N_2 \oplus \M_2 \oplus \M_1 $ to $ \N_1 \oplus \M_1 \oplus \M_2$ sending a generic element $n_2 \oplus m_2 \oplus m_1$ to $Vn_2 \oplus m_1 \oplus m_2.$ Taking $Vn^\prime _2 =n^\prime _1$, $\mathbf{t}(z)^* m^\prime _2=m^\prime _1$, $n^\prime _1 \oplus m^\prime _1= q_z \oplus u$ and $n^\prime _2\oplus m^\prime _2=  \lamZ^* q_z \oplus y$ we see that $\Q$ sends  $$\lamZ^* q_z \oplus  \mathbf{t}(z)^* m^\prime _2\oplus y\,\, \text{to}\,\, q_z \oplus m^\prime _2\oplus u.$$
Hence we are allowed to write
\begin{align*}
\bordermatrix{ & \mathscr{L}_1 & \M_1 \oplus \Y  \cr
	\qquad \mathscr{L}_1 & \Q_{11} & \Q_{12} \cr
	\M_2 \oplus \U & \Q_{12} & \Q_{22} } \begin{pmatrix} \lamZ^* q_z \\ \mathbf{t}(z)^* m^\prime _2\oplus y \end{pmatrix} & = \begin{pmatrix} q_z \\ m^\prime _2\oplus u \end{pmatrix}.
\end{align*}
Hence
$$
\Q_{11}(\lamZ^* q_z) + \Q_{12} (\mathbf{t}(z)^* m^\prime _2\oplus y) = q_z$$
$$\Q_{21} (\lamZ^* q_z) + \Q_{22} (\mathbf{t}(z)^* m^\prime _2\oplus y) = m^\prime _2\oplus u.
$$
Eliminating $q_z$ we obtain
$$(\Q_{22} + \Q_{21} \lamZ^*(I_{\LL_1} - \Q_{11} \lamZ^* )^{-1} \Q_{12})(\mathbf{t}(z)^* m^\prime _2\oplus y)= m^\prime _2 \oplus u.$$
Now recall that from (\ref{equality of lambda and mu}) we have
\begin{align*}
\lambda (\delta)|_{\LL_1} = \mu(\delta)|_{\LL_1}\,\, \text{for all}\,\, \delta \in \CB.
\end{align*}
So we have
$$\Big(\Q_{22} + \Q_{21} \mu (E(z))^*(I_{\LL_1} - \Q_{11} \mu (E(z))^* )^{-1} \Q_{12}\Big)(\mathbf{t}(z)^* m^\prime _2\oplus y)= m^\prime _2 \oplus u.$$
Recalling the $G$ from (\ref{defn_G}) of Section \ref{Affiliated}, we see that the last equation is precisely
$$G(z)^* (\mathbf{t}(z)^* m^\prime _2\oplus y)=m^\prime _2 \oplus u.$$
Using the decomposition of $G(z)^*$ with respect to its domain and range we get
\begin{align*}
\bordermatrix{ & \M_1 &  \Y  \cr
	\M_2 & G_{11}(z)^* & G_{21}(z)^* \cr
	\U & G_{12}(z)^* &G_{22}(z)^* } \begin{pmatrix}\mathbf{t}(z)^* m^\prime _2  \\ y \end{pmatrix} =  \begin{pmatrix} m^\prime _2  \\ u \end{pmatrix}.
\end{align*}
From this we obtain two equations. Eliminating $m_2$ from those gives us
$$u= G_{22}(z)^*y + G_{12}(z)^*\mathbf{t}(z)^*(I_{\M_2}-G_{11} (z)^* \mathbf{t}(z)^* )^{-1}G_{21} (z)^* y.$$
Since we know from Section \ref{Affiliated} that $G_{11} (z)<1$ for each $z\in \Omega$, the right hand side is well defined. As $u=f(z)^* y$ and, $y$ and $z$ are arbitrary, we have
$$f(z) = G_{22}(z) + G_{21}(z) (I_{\M_2} - \mathbf{t}(z) G_{11}(z))^{-1} \mathbf{t}(z) G_{12}(z),\,\, \text{for all} \,\, z\in \Omega.$$
This completes the proof.

\section{Examples}\label{Examples}

The {\em Schur class} $H_1 ^\infty (\Omega, B(\U,\Y))$ of a domain $\Omega$ is the closed unit ball (in the supremum norm) of the algebra of all bounded holomorphic functions on $\Omega$ taking values in $B(\U,\Y)$.

\begin{theorem}
Suppose $\Omega$ stands for the bidisc or the symmetrized bidisc or the annulus. We consider two Hilbert spaces $\U$ and $\Y$. When $\Omega$ is the annulus, we take $\U =\Y =\mathbb{C}$. 	 Suppose that $f \in H_1 ^\infty (\Omega, B(\U,\Y))$ is a solution of this interpolation problem $z_i \mapsto B_i$ and $f$ is affiliated with a completely positive kernel
	\begin{align*}
	\Gamma: \{z_1,z_2,\ldots, z_n\} \times  \{z_1,z_2,\ldots, z_n\}\rightarrow B(\mathcal{C}_b(\Psi), B(\Y)).
	\end{align*}
	Then with $\M_1$, $ \M_2$ and $G$ as in Section \ref{Affiliated}, we have that writing $G$ as
	\[
	G(z) =
	\bordermatrix{ & \mathscr{M}_2 & \U \cr
		\mathscr{M}_1 & G_{11}(z) & G_{12}(z) \cr
		\Y & G_{21}(z) & G_{22}(z)}, \qquad
	\]
	one has $$f(z) = \Big(G_{22} + G_{21}(I_{\M_2} -\mathbf{t}\  G_{11})^{-1} \mathbf{t}\  G_{12}\Big)(z)$$ for some $\mathbf{t} \in H_1 ^\infty (\Omega, B(\M_1, \M_2))$ and for all $ z \in \Omega.$
\end{theorem}

\textbf{Proof} \\

 In each of these examples, there exists a certain collection of holomorphic test functions, say $\Psi$, which satisfies the fact that there is a point $w_0$ in the domain where all test functions vanish. We do not get into the details of writing down the test functions explicitly for the sake of brevity. While for the bidisc the collection consists of just two test functions - the co-ordinate functions $z_1$ and $z_2$, for the symmetrized bidisc and the annulus, they are uncountable in number. See \cite{B-S} for the symmetrized bidisc and \cite{D-M} for the annulus.  Now we apply our Main Theorem. The crucial fact which clinches the issue is that in each of the above domains, for the above mentioned test functions, the $\Psi$-Schur-Agler Class coincides with the Schur class $H_1 ^\infty (\Omega, B(\U,\Y))$. \qed

\vspace{0.1in} \noindent\textbf{Acknowledgement:}
The first named author's research is supported by the University Grants Commission Centre for Advanced Studies. The third named author's research is supported by a Post Doctoral Fellowship at Indian Institute of Technology, Bombay. All authors thank the referee because the referee's inputs led to a substantial increase in the quality of the paper.

\end{document}